\documentclass[11pt]{article}
\usepackage{amsfonts}
\usepackage{mathrsfs}
\usepackage{amsmath}
\usepackage{amssymb}
\usepackage{graphicx}
\usepackage{epsfig}
\usepackage{epic}
\usepackage{cite}
\usepackage{setspace}
\usepackage{amsthm,latexsym}
\usepackage{float}

\renewcommand{\paragraph}{\roman{paragraph}}
\newtheorem{theorem}{\scshape \mdseries \bf Theorem}[section]
\newtheorem{lemma}[theorem]{\scshape \mdseries  \bf Lemma}
\newtheorem{remark}[theorem]{\scshape \mdseries  \bf Remark}
\newtheorem{corol}[theorem]{\scshape \mdseries  \bf Corollary}

\begin{document}

\title{\sf Characterization of graphs with some normalized Laplacian eigenvalue of multiplicity $n-3$}
\author{ \ Fenglei Tian$^{1,}$\thanks{Corresponding author. E-mail address: tflqsd@126.com.  Supported by '' the Natural Science Foundation of Shandong Province (No. ZR2019BA016)
''. }\ ,\
Dein Wong$^2$\\
~~\\
\noindent{\small\it 1.\ School of Management, Qufu Normal University, Rizhao, China.}\\
\noindent{\small\it 2.\ School of Mathematics, China University of Mining and Technology, Xuzhou, China.}
   }
\date{}
\maketitle
\noindent {\bf Abstract:} \ Graphs with few distinct eigenvalues have been investigated extensively. In this paper, we focus on another relevant topic: characterizing graphs with some eigenvalue of large multiplicity. Specifically, the normalized Laplacian matrix of a graph is considered here. Let $\rho_{n-1}(G)$ and $\nu(G)$ be the second least normalized Laplacian eigenvalue and the independence number of a graph $G$, respectively.
As the main conclusions, two families of $n$-vertex connected graphs with some normalized Laplacian eigenvalue of multiplicity $n-3$ are determined: graphs with $\rho_{n-1}(G)=-1$ and graphs with $\rho_{n-1}(G)\neq -1$ and $\nu(G)\neq 2$. Moreover, it is proved that these graphs are determined by their spectrum.

\vskip 2 mm
\noindent{\bf Keywords:}\ Normalized Laplacian matrix; Normalized Laplacian eigenvalues; Eigenvalue multiplicity
\vskip 1.5 mm
\noindent{\bf AMS classification:}\ \ 05C50

\section{Introduction}

\quad {Investigating graphs with few distinct eigenvalues has attracted much attention, since it was proposed. Several matrices associated with graphs have been considered, such as the adjacency matrix \cite{Doob,vanDam1,vanDam2,vanDam3,vanDam4,Muzychuk,Rowlinson,Cheng1,Cheng2,Huang1}, the Laplacian matrix \cite{vanDam5,wangyi,Mohammadian}, the signless Laplacian matrix \cite{Ayoobi}, the universal adjacency matrix \cite{Haemers} and the normalized Laplacian matrix \cite{vanDam6,Braga,Huang2}.
One of the reasons to study such graphs is that they seem to have certain kind of regularity \cite{Brouwer}. Moreover, another motivation for considering graphs with few distinct eigenvalues is that most of those graphs are not determined by their spectrum \cite{vanDam6}. Hence, it is related to the question of which graphs are determined by their spectrum \cite{vanDam7}. However, for any kind of matrices associated with graphs, it is not easy to give a complete characterization for such graphs. Therefore, searching more families of graphs with few distinct eigenvalues is of interest.

The normalized Laplacian spectrum of graphs has been studied intensively, as it reveals some structural properties and some relevant dynamical aspects (such as random walk) of graphs \cite{Chung}. But the results on graphs with few distinct normalized Laplacian eigenvalues are restricted. van Dam and Omidi \cite{vanDam6} first gave a combinatorial characterization of graphs with three distinct normalized Laplacian eigenvalues and constructed some special families of such graphs. Braga $et\ al.$ \cite{Braga} characterized trees with 4 or 5 distinct normalized Laplacian eigenvalues. Huang $et\ al.$ determined all connected graphs having three distinct normalized Laplacian eigenvalues with one equal to 1, and determined other classes of graphs with three or four distinct normalized Laplacian eigenvalues.

In some sense, investigating the graphs with some eigenvalue of large multiplicity is in connection with characterizing the graphs with few distinct eigenvalues. For instance, assume that $G$ is a connected graph of order $n>3$, then $G$ has a normalized Laplacian eigenvalue with multiplicity $n-1$ if and only if $G$ has two distinct normalized Laplacian eigenvalues (note that $0$ as a normalized Laplacian eigenvalue is simple for connected graphs); $G$ has a normalized Laplacian eigenvalue with multiplicity $n-2$ if and only if $G$ has three distinct normalized Laplacian eigenvalues and two of them are simple.
In fact, van Dam and Omidi [16, Proposition 8] has determined the graphs with some normalized Laplacian eigenvalue of multiplicity $n-2$. In this paper, we further focus on the connected graphs with some normalized Laplacian eigenvalue of multiplicity $n-3$ as an extension of the result in [16], obtaining the following conclusion.

Denote by $\mathcal{G}(n, n-3)$ the set of all $n$-vertex ($n\geq 5$) connected graphs with some normalized Laplacian eigenvalue of multiplicity $n-3$. It is well-known that the least normalized Laplacian eigenvalue of a connected graph is $0$ with multiplicity $1$ (see \cite{Chung}). Then let the eigenvalues of the normalized Laplacian matrix of a connected graph $G$ be
$$\rho_1(G)\geq \rho_2(G)\geq \cdots \geq \rho_{n-1}(G)> \rho_{n}(G)=0.$$
The independence number of $G$ is denoted by $\nu(G)$.

\begin{theorem}\ Let $G$ be a connected graph of order $n\geq 5$. Then \end{theorem}
\begin{spacing}{1}
\begin{enumerate}
\item[(i)]  \textit{$G\in \mathcal{G}(n, n-3)$ with  $\rho_{n-1}(G)=1$ if and only if $G$ is a complete tripartite graph $K_{a,b,c}$ or $K_n-e$, where $K_n-e$ is the graph obtained from the complete graph $K_n$ by removing an edge.}
\item[(ii)]  \textit{$G\in \mathcal{G}(n, n-3)$ with  $\rho_{n-1}(G)\neq 1$ and $\nu(G)\neq 2$ if and only if $G\in \{G_1, G_2, G_3\}$ (see Fig. 1).}
\end{enumerate}
\end{spacing}

\begin{figure}[htbp]
  \centering
  \setlength{\abovecaptionskip}{0cm} 
  \setlength{\belowcaptionskip}{0pt}
  \includegraphics[width=5 in]{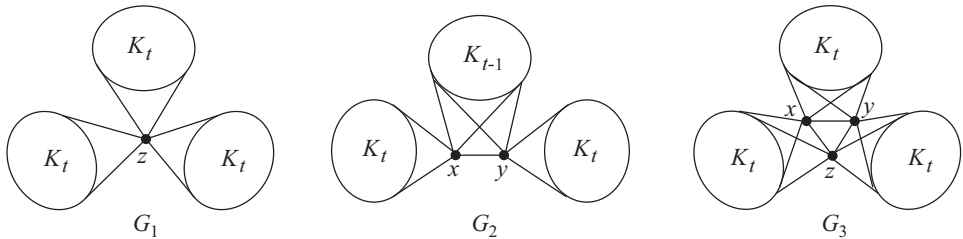}
  \caption{The graphs $G_1$, $G_2$ and $G_3$.}
\end{figure}

Before showing the proof of Theorem 1.1, we first introduce some notations and lemmas in the following Section 2.

\section{ Preliminaries }

\quad
Throughout, all graphs considered here are connected and simple. Let $G=(V(G),E(G))$ be a graph with  vertex set $V(G)$ and edge set $E(G)$. The set of all the neighbors of a vertex $u$ is denoted by $N_G(u)$, and $d_u=|N_G(u)|$ is called the degree of $u$. For a subset $S$ of $V(G)$, $S$ is called a set of twin points if $N_G(u)=N_G(v)$ for any $u, v\in S$. The notation $u\thicksim v$ means that $u$ is adjacent to $v$. If any two vertices of a subset of $V(G)$ are nonadjacent, then the subset is called an independent set of $G$. The independence number $\nu(G)$ of $G$ is the cardinality of the maximum independent set. The rank of a matrix $M$ is denoted by $r(M)$.

Let $A(G)$ and $L(G)=D(G)-A(G)$ be the adjacency matrix and the Laplacian matrix of a graph $G$, respectively. Then the normalized Laplacian matrix $\mathcal{L}(G)=[l_{uv}]$ of graph $G$ is defined as $$\mathcal{L}(G)=D^{-1/2}(G)L(G)D^{-1/2}(G)=I-D^{-1/2}(G)A(G)D^{-1/2}(G),$$
where
$$l_{uv}=\begin{cases}
 \ 1, \ \ \ \ \ \ \ \ \ \ \ \ \ \text{if \ $u=v$};\\
 -1/\sqrt{d_ud_v}, \ \ \text{if $u\thicksim v$};\\
 \ 0, \ \ \ \ \ \ \ \ \ \ \ \ \ \text{otherwise}.
 \end{cases}$$
For brevity, the normalized Laplacian eigenvalue is written as $\mathcal{L}$-eigenvalue.
In what follows, some known results are listed.


\begin{lemma}\label{Interlacinglemma}{\rm (Interlacing Theorem, \cite{Brouwer})}\  Let $A$ be a real and symmetric matrix of order $n$ and $M$ be a principal submatrix of $A$ with order $s(\leq n)$. Then
\begin{eqnarray}
 \nonumber \lambda_{i+n-s}(A)\leq \lambda_i(M)\leq \lambda_i(A),\ \  1\leq i\leq s,
 \end{eqnarray}
where $\lambda_i$ denotes the $i$-th largest eigenvalue.
\end{lemma}

Let $H$ be a real symmetric matrix of order $n$ whose columns and rows are indexed by $X=\{1, 2, \cdots, n\}$. Denote by $\{X_1, X_2, \cdots, X_t\}$ a partition of $X$. An $n$-dimensional column vector whose components indexed by $X_i$ are ones and all others are zeros is called the characteristic vector of $X_i$. The $n\times t$ matrix $S$ whose $i$-th column is the characteristic vector of $X_i$ is called the characteristic matrix. Let the block form of $H$ with respect to the partition of $X$ be
$$H=
\left(
  \begin{array}{ccc}
    H_{11} & \cdots & H_{1t} \\
    \vdots & \ddots & \vdots \\
    H_{t1} & \cdots & H_{tt} \\
  \end{array}
\right),$$
where $H_{ji}$ is the transpose of $H_{ij}$.
Let $q_{ij}$ be the average row sum of $H_{ij}$, then the matrix $Q=(q_{ij})$ is called the quotient matrix of $H$. Moreover, if the row sum of $H_{ij}$ is constant, then the partition of $X$ is called equitable (see \cite{Brouwer}).

\begin{lemma}\label{quotientmatrixlemma}{\rm \cite{Brouwer}} \ Suppose that $H$ is a real symmetric matrix with an equitable partition. Let $S$ and $Q$ be the corresponding characteristic matrix and quotient matrix of $H$, respectively. Then, each eigenvalue of $Q$ is an eigenvalue of $H$. Moreover, if $\alpha$ is an eigenvector of $Q$ for an eigenvalue $\lambda$, then $S\alpha$ is an eigenvector of $H$ for the eigenvalue $\lambda$.
\end{lemma}


\begin{lemma}\label{twinpointslemma}{\rm \cite{Huang2,Das}}\  Let $G$ be a graph with order $n$. Denote by $\{v_1,\ldots,v_p \}$ a set of twin points of $G$, then $1$ is an $\mathcal{L}$-eigenvalue of $G$ with multiplicity at least $p-1$.
\end{lemma}

\begin{lemma}\label{cliquelemma}{\rm\cite{Huang2}}\  Let $G$ be a graph with $n$ vertices. Let $K=\{v_1, \ldots, v_q \}$ be a clique in $G$ such that $N_G(v_i)-K=N_G(v_j)-K$ $(1 \leq i, j\leq q)$, then $1+\frac{1}{d_{v_i}}$ is an $\mathcal{L}$-eigenvalue of $G$ with multiplicity at least $q-1$. Further, the $q-1$ eigenvectors for $1+\frac{1}{d_{v_i}}$ can be written as
\vspace{-0.6cm}
$$\alpha_1=
\bordermatrix{
    &   \cr
    & 1  \cr
2   & -1  \cr
    & 0  \cr
    & \vdots \cr
    & \vdots \cr
    & \vdots \cr
    & \vdots \cr
    & 0 \cr
},\ \
\alpha_2=
\bordermatrix{
    &   \cr
    & 1  \cr
    & 0  \cr
3   & -1  \cr
    & 0  \cr
    & \vdots \cr
    & \vdots \cr
    & \vdots \cr
    & 0 \cr
},
\cdots\cdots,\ \
\alpha_{q-1}=
\bordermatrix{
    &   \cr
    & 1  \cr
    & 0  \cr
    & \vdots \cr
    & 0  \cr
q   & -1  \cr
    & 0  \cr
    & \vdots \cr
    & 0 \cr}
.$$
\end{lemma}



\begin{lemma}\label{secondlargesteigenvalue}{\rm\cite{Guo}}\ Let $G$ be a graph with $n(\geq 3)$ vertices. Then $\rho_2(G)\geq 1$ and equality holds if and only if $G$ is a complete bipartite graph.
\end{lemma}

\begin{lemma}\label{secondleasteigenvalue}{\rm\cite{Guo}}\ Let $G$ be a graph with $n$ vertices, which is not a complete graph. Then $\rho_{n-1}(G)\leq 1$ and equality holds if and only if $G$ is a complete multipartite graph.
\end{lemma}

\begin{lemma}\label{less3}\ Let $G\in\mathcal{G}(n, n-3)$ and $\theta$ be the $\mathcal{L}$-eigenvalue of $G$ with multiplicity $n-3$. If $\rho_{n-1}(G)\neq 1$, then $\theta\neq 1$ and $\nu(G)\leq 3$.
\end{lemma}

\noindent
{\bf Proof.} \  It is obvious that the rank $r(\mathcal{L}(G)-\theta I)$ of $\mathcal{L}(G)-\theta I$ is 3, where $I$ is the identity matrix. Suppose on the contrary that $\theta= 1$, then $\rho_2(G)=\theta=1$ (noting that $\rho_{n-1}(G)\neq 1$ and $\rho_{n}(G)= 0$). Thus $G$ is a complete bipartite graph from Lemma \ref{secondlargesteigenvalue}.
However, from $\rho_{n-1}(G)\neq 1$ and Lemma \ref{secondleasteigenvalue}, $G$ is not a complete multipartite graph, a contradiction. Hence, $\theta\neq 1$ holds.
Now, assume that $\nu(G)> 3$ and $\{u, v, w, z\}$ is an independent set of $G$. Let $M$ be the principal submatrix of $\mathcal{L}(G)-\theta I$ indexed by $\{u, v, w, z\}$. Then $r(\mathcal{L}(G)-\theta I)\geq r(M)= 4$, a contradiction. \hfill$\square$

\begin{lemma}\label{commonvertex}\ Let $G\in\mathcal{G}(n, n-3)$ with $\nu(G)=3$ and $\rho_{n-1}(G)\neq 1$. Suppose that $\{u, v, w\}$ is an independent set of $G$. Then the following assertions hold.\end{lemma}
\begin{spacing}{1}
\begin{enumerate}
\item[(i)]  \textit{If there exists a vertex, say $x$, such that $x\thicksim v$, $x\thicksim w$ and $x\nsim u$. Then there exists no other vertex distinct with $x$, which is adjacent to $v$ and $w$.}
\item[(ii)]  \textit{If $z$ is a vertex adjacent to exactly one of $\{u, v, w\}$, say $u$, then $N_G(z)-u=N_G(u)-z$.}
\item[(iii)]  \textit{Any of $\{u, v, w\}$ must have a common vertex with at least one of the remaining two of $\{u, v, w\}$.}
\item[(iv)]  \textit{There exists at most one vertex adjacent to each of $\{u, v, w\}$.}
\item[(v)]  \textit{If there exist two vertices, say $x$ and $y$, such that $x\thicksim u$, $x\thicksim v$, $x\nsim w$ and $y\nsim u$, $y\thicksim v$, $y\thicksim w$. Then $x\thicksim y$, and $d^2_v=d_ud_x-d_ud_v=d_yd_w-d_wd_v$.}
\end{enumerate}
\end{spacing}

\noindent
{\bf Proof.} \  Denote by $R_u$ the row of $\mathcal{L}(G)-\theta I$ indexed by the vertex $u$. Then one can easily see that the rows $\{R_u, R_v, R_w\}$ of $\mathcal{L}(G)$ indexed by $\{u, v, w\}$ are linearly independent. Thus any other row of $\mathcal{L}(G)-\theta I$ can be written as a linear combination of $\{R_u, R_v, R_w\}$.

\vskip 2mm
First, we show the proof of assertion (i).
From the above discussion, there exist three real numbers $\{a, b, c\}$ such that
\begin{equation}\label{e1}
R_x=aR_u+bR_v+cR_w.
\end{equation}
Suppose for a contradiction that there exists another vertex distinct with $x$, say $y$, such that $y\thicksim v$ and $y\thicksim w$.
Let $M$ be the principal submatrix of $\mathcal{L}(G)-\theta I$ indexed by the vertices $\{u, v, w, x, y\}$, then
$$M=\left(
       \begin{array}{ccccc}
         1-\theta & 0 & 0 & 0 & \delta_1 \\
         0 & 1-\theta & 0 & \frac{-1}{\sqrt{d_vd_x}} & \frac{-1}{\sqrt{d_vd_y}} \\
         0 & 0 &  1-\theta & \frac{-1}{\sqrt{d_wd_x}} &\frac{-1}{\sqrt{d_wd_y}} \\
         0 & \frac{-1}{\sqrt{d_xd_v}} &  \frac{-1}{\sqrt{d_xd_w}} & 1-\theta & \delta_2 \\
         \delta_1 & \frac{-1}{\sqrt{d_yd_v}} & \frac{-1}{\sqrt{d_yd_w}} & \delta_2 & 1-\theta \\
       \end{array}
     \right),$$
where $\delta_1=0$ or $\frac{-1}{\sqrt{d_ud_y}}$ with respect to $u\nsim y$ or $u\thicksim y$ and $\delta_2=0$ or $\frac{-1}{\sqrt{d_xd_y}}$ with respect to $x\nsim y$ or $x\thicksim  y$.

Since $\theta\neq 1$ from Lemma \ref{less3}, then applying Eq. (\ref{e1}) to the first three columns of $M$, we obtain that
$$a=0,\ ~\ b= \frac{-1}{(1-\theta)\sqrt{d_xd_v}},\ ~\ c= \frac{-1}{(1-\theta)\sqrt{d_xd_w}}.$$
Hence,
\begin{equation}\label{e2}
R_x=\frac{-1}{(1-\theta)\sqrt{d_xd_v}}R_v+\frac{-1}{(1-\theta)\sqrt{d_xd_w}}R_w.
\end{equation}
From Eq. (\ref{e2}), if $z$ is a vertex of $G$, then $z\thicksim x$ if and only if $z\thicksim v$ or $z\thicksim w$, which implies that $$d_x>d_v.$$
As a result, $y\thicksim x$ and $\delta_2=\frac{-1}{\sqrt{d_xd_y}}$.
Now, applying Eq. (\ref{e2}) to the fourth and the fifth columns of $M$, it follows that
\begin{equation*}
\begin{cases}
 \frac{1}{(1-\theta)d_vd_x}+\frac{1}{(1-\theta)d_wd_x}=1-\theta\\
 \frac{1}{(1-\theta)d_v}+\frac{1}{(1-\theta)d_w}=-1,
 \end{cases}
\end{equation*}
which yields
$$\frac{1}{d_x}=\frac{1}{d_v}+\frac{1}{d_w},$$
contradicting with $d_x>d_v$.
Consequently, there exists no vertex $y$ such that $y\thicksim v$ and $y\thicksim w$.

\vskip 2mm
For assertion (ii), similar as Eq. (\ref{e1}), let $R_z=aR_u+bR_v+cR_w$. Then from the principal submatrix of $\mathcal{L}(G)-\theta I$ indexed by $\{u, v, w, z\}$, we easily obtain $b=c=0$ and $a\neq 0$, then the result is clear.

\vskip 2mm
For assertion (iii), without loss of generality, suppose that $u$ has no common vertex with both $v$ and $w$. Recalling that $\nu(G)=3$, any vertex out of $\{u, v, w\}$ must be adjacent to at least one of $\{u, v, w\}$. From assertion (ii), the vertices adjacent to $u$ also have no common vertex with $v$ and $w$. Then $G$ is not connected, a contradiction.

\vskip 2mm
Next, for assertion (iv), assume that there exist two vertices, say $s$ and $t$, adjacent to each of $\{u, v, w\}$. Then from the assertion (i), there exists no vertex adjacent to exactly two of $\{u, v, w\}$.

Let $M$ be the principal submatrix of $\mathcal{L}(G)-\theta I$ indexed by the vertices $\{u, v, w, s, t\}$, then
$$M=\left(
       \begin{array}{ccccc}
         1-\theta & 0 & 0 & \frac{-1}{\sqrt{d_ud_s}} & \frac{-1}{\sqrt{d_ud_t}} \\
         0 & 1-\theta & 0 & \frac{-1}{\sqrt{d_vd_s}} & \frac{-1}{\sqrt{d_vd_t}} \\
         0 & 0 & 1-\theta & \frac{-1}{\sqrt{d_wd_s}} & \frac{-1}{\sqrt{d_wd_t}} \\
         \frac{-1}{\sqrt{d_sd_u}} & \frac{-1}{\sqrt{d_sd_v}} &  \frac{-1}{\sqrt{d_sd_w}} & 1-\theta & \delta \\
         \frac{-1}{\sqrt{d_td_u}} & \frac{-1}{\sqrt{d_td_v}} & \frac{-1}{\sqrt{d_td_w}} & \delta & 1-\theta \\
       \end{array}
     \right),$$
where $\delta=0$ or $\frac{-1}{\sqrt{d_sd_t}}$ according to $s\nsim t$ or $s\thicksim t$, respectively. But we claim that $s\thicksim t$, i.e., $\delta=\frac{-1}{\sqrt{d_sd_t}}$. The following is the reason. Let $R_s=aR_u+bR_v+cR_w$, then it follows from the first three columns of $M$ and $\theta\neq 1$ that
\begin{equation}\label{e3}
\begin{cases}
a=\frac{-1}{(1-\theta)\sqrt{d_sd_u}},\\ b=\frac{-1}{(1-\theta)\sqrt{d_sd_v}},\\ c=\frac{-1}{(1-\theta)\sqrt{d_sd_w}}.
\end{cases}
\end{equation}
Further, by the fifth column of $M$,
\begin{equation}\label{e4}
a\frac{-1}{\sqrt{d_ud_t}}+b\frac{-1}{\sqrt{d_vd_t}}+c\frac{-1}{\sqrt{d_wd_t}}=\delta,
\end{equation}
which, together with Eq. (\ref{e3}), implies that $\delta\neq 0$. Thus $s\thicksim t$, i.e., $\delta=\frac{-1}{\sqrt{d_sd_t}}$.
Now, simplifying Eq. (\ref{e4}), we obtain that
\begin{equation}\label{e5}
\frac{1}{d_u}+\frac{1}{d_v}+\frac{1}{d_w}=-(1-\theta).
\end{equation}
From the fourth column of $M$ and Eq. (\ref{e3}), we derive that
\begin{equation}\label{e6}
\frac{1}{d_s}(\frac{1}{d_u}+\frac{1}{d_v}+\frac{1}{d_w})=(1-\theta)^2.
\end{equation}
Combining Eqs. (\ref{e5}) and (\ref{e6}), it follows that
\begin{equation}\label{q1}
\frac{1}{d_s}=\frac{1}{d_u}+\frac{1}{d_v}+\frac{1}{d_w}.
\end{equation}

Recalling that there is no vertex adjacent to exactly two of $\{u, v, w\}$ and any two vertices adjacent to each of $\{u, v, w\}$ must be adjacent, from assertion (ii) one can obtain
$$d_s>d_u,$$
contradicting with Eq. (\ref{q1}).
Consequently, the assertion (iv) holds.   \hfill$\square$

\vskip 2mm
At last, we prove assertion (v).
Suppose that $x, y$ are the vertices as stated in the condition of assertion (v). Let $M$ be the principal submatrix of $\mathcal{L}(G)-\theta I$ indexed by the vertices $\{u, v, w, x, y\}$. If $x\nsim y$, then the vertices $\{u, v, w, x, y\}$ induce a path $P_5$. Then it is easy to see that $r(\mathcal{L}(G)-\theta I)\geq r(M)\geq 4$, contradicting with $r(\mathcal{L}(G)-\theta I)=3$. Hence, we say that $x\thicksim y$, and the principal submatrix $M$ of $\mathcal{L}(G)-\theta I$ can be written as
$$M=\left(
       \begin{array}{ccccc}
         1-\theta & 0 & 0 &\frac{-1}{\sqrt{d_ud_x}} & 0 \\
         0 & 1-\theta & 0 & \frac{-1}{\sqrt{d_vd_x}} & \frac{-1}{\sqrt{d_vd_y}} \\
         0 & 0 & 1-\theta & 0 &\frac{-1}{\sqrt{d_wd_y}} \\
         \frac{-1}{\sqrt{d_xd_u}} & \frac{-1}{\sqrt{d_xd_v}} & 0 & 1-\theta & \frac{-1}{\sqrt{d_xd_y}} \\
         0 & \frac{-1}{\sqrt{d_yd_v}} & \frac{-1}{\sqrt{d_yd_w}} & \frac{-1}{\sqrt{d_yd_x}} & 1-\theta \\
       \end{array}
     \right).$$
Similar as Eq. (\ref{e1}), we set
\begin{equation}\label{e7}
R_{x}=aR_{u}+bR_{v}+cR_{w}.
\end{equation}
Applying Eq. (\ref{e7}) to the third column of $M$, we get $c=0$. Further, applying Eq. (\ref{e7}) to the remaining columns of $M$, it is obtained that
\begin{equation}\label{e8}
\begin{cases}
 a(1-\theta)=\frac{-1}{\sqrt{d_xd_u}},\\
 b(1-\theta)=\frac{-1}{\sqrt{d_xd_v}},\\
 \frac{-a}{\sqrt{d_ud_x}}-\frac{b}{\sqrt{d_vd_x}}=1-\theta,\\
 \frac{-b}{\sqrt{d_vd_y}}=\frac{-1}{\sqrt{d_xd_y}}.
 \end{cases}
\end{equation}
The second and the fourth equations of (\ref{e8}) tell us that $1-\theta=-\frac{1}{d_v},$
and the first three equations of (\ref{e8}) imply that $$(1-\theta)^2=\frac{1}{d_x}(\frac{1}{d_u}+\frac{1}{d_v}).$$
Then combining the two equations above, we get $d_v^2=d_ud_x-d_ud_v$. Analogously, one can derive that $d_v^2=d_yd_w-d_wd_v$.

Consequently, all the proofs are completed. \hfill$\square$

\begin{lemma}\label{spectrum}\ Let $G_1, G_2$ and $G_3$ be the graphs in Fig. 1. Then their spectra (eigenvalues with multiplicity) are respectively
\begin{equation*}
\begin{cases}
\{0^1,\ {(\frac{3}{n-1})}^2,\ {(\frac{n+2}{n-1})}^{n-3}\},\\
\{0^1,\ {(\frac{3}{2(n-1)})}^1,\ {(\frac{9}{2(n-1)})}^1,\  {(\frac{n+2}{n-1})}^{n-3}\},\\
\{0^1,\ {(\frac{9}{2n})}^2,\ {(\frac{n+3}{n})}^{n-3}\}.
\end{cases}
\end{equation*}
\end{lemma}

\noindent
{\bf Proof.} \ For graph $G_1$, divide its vertex set $V(G_1)$ into four parts
$$V(G_1)=V(K_t)\cup V(K_t)\cup V(K_t)\cup \{z\}.$$
Accordingly, $\mathcal{L}(G_1)$ has an equitable partition. Let $Q_1$ be the quotient matrix of $\mathcal{L}(G_1)$, then
$$Q_1=\left(
  \begin{array}{cccc}
    1-\frac{t-1}{t} & 0 & 0 & -\frac{1}{\sqrt{t(3t)}} \\
    0 & 1-\frac{t-1}{t} & 0 & -\frac{1}{\sqrt{t(3t)}} \\
    0 & 0 & 1-\frac{t-1}{t} & -\frac{1}{\sqrt{t(3t)}} \\
    -\frac{t}{\sqrt{t(3t)}} & -\frac{t}{\sqrt{t(3t)}} & -\frac{t}{\sqrt{t(3t)}} & 1 \\
  \end{array}
\right).$$
By direct calculation, the eigenvalues of $Q_1$ are
$$\{0,\ \frac{t+1}{t},\ \frac{1}{t},\ \frac{1}{t}\},$$
which are also the eigenvalues of $\mathcal{L}(G_1)$ from Lemma \ref{quotientmatrixlemma}. Furthermore, $\frac{t+1}{t}$ is an $\mathcal{L}$-eigenvalue of $G_1$ with multiplicity at least $3t-3$ by Lemma \ref{cliquelemma}. Noting that $3t+1=n$, it follows from the trace of $\mathcal{L}(G_1)$ that the last $\mathcal{L}$-eigenvalue of $G_1$ is $$n-(3t-3)\cdot \frac{t+1}{t}-2\cdot \frac{1}{t}=
n-(n-4)\cdot \frac{n+2}{n-1}-2\cdot \frac{3}{n-1}=\frac{n+2}{n-1},$$
which implies that the multiplicity of the $\mathcal{L}$-eigenvalue $\frac{n+2}{n-1},\ (i.e., \frac{t+1}{t})$ is $n-3$. Thus the spectrum of $G_1$ is $\{0^1,\ {(\frac{3}{n-1})}^2,\ {(\frac{n+2}{n-1})}^{n-3}\}$.

For graph $G_2$, $\mathcal{L}(G_2)$ has an equitable partition according to the vertex partition
$$V(G_2)=V(K_t)\cup V(K_{t-1})\cup V(K_t)\cup \{x\}\cup \{y\}.$$
Denote the corresponding characteristic matrix and quotient matrix of $\mathcal{L}(G_2)$ by $S_2$ and $Q_2$ respectively, then
\vspace{-0.6cm}
$$S_2=\ \bordermatrix{
    & & & & &  \cr
1      & 1        & 0       & 0       & 0      & 0        \cr
       & \vdots   & \vdots  & \vdots  & \vdots & \vdots   \cr
t      & 1        & 0       & 0       & 0      & 0        \cr
(t+1)  & 0        & 1       & 0       & 0      & 0        \cr
       & \vdots   & \vdots  & \vdots  & \vdots & \vdots  \cr
(2t-1) & 0        & 1       & 0       & 0      & 0       \cr
2t     & 0        & 0       & 1       & 0      & 0      \cr
       & \vdots   & \vdots  & \vdots  & \vdots & \vdots  \cr
(3t-1) & 0        & 0       & 1       & 0      & 0      \cr
3t     & 0        & 0       & 0       & 1      & 0     \cr
(3t+1) & 0        & 0       & 0       & 0      & 1     \cr}
,\ \
Q_2=
\left(
  \begin{array}{ccccc}
    \frac{1}{t} & 0 & 0 & -\frac{1}{t\sqrt{2}} & 0 \\
    0 & \frac{2}{t} & 0 & -\frac{1}{t\sqrt{2}} & -\frac{1}{t\sqrt{2}} \\
    0 & 0 & \frac{1}{t} & 0 & -\frac{1}{t\sqrt{2}} \\
    -\frac{1}{\sqrt{2}} & -\frac{t-1}{t\sqrt{2}} & 0 & 1 & -\frac{1}{2t} \\
    0 & -\frac{t-1}{t\sqrt{2}} & -\frac{1}{\sqrt{2}} & -\frac{1}{2t} & 1 \\
  \end{array}
\right).$$
By direct calculation, the eigenvalues of $Q_2$ are
$$\{0,\ \frac{1}{2t},\ \frac{3}{2t},\ \frac{1+t}{t},\  \frac{1+t}{t}\},$$
and there are two linearly independent eigenvectors, denoted by $\xi_1$ and $\xi_2$, for the eigenvalue $\frac{1+t}{t}$.
From Lemma \ref{quotientmatrixlemma}, all the eigenvalues of $Q_2$  are also the eigenvalues of $\mathcal{L}(G_2)$, and $S_2\xi_1,\ S_2\xi_2$ are two linearly independent eigenvectors for the  $\mathcal{L}$-eigenvalue $\frac{1+t}{t}$.
Further, by Lemma \ref{cliquelemma}, $\frac{1+t}{t}$ is an $\mathcal{L}$-eigenvalue of $G_2$ with multiplicity at least $3t-4$ and the corresponding $3t-4$ eigenvectors for $\frac{1+t}{t}$ can be easily obtained similar as stated in Lemma \ref{cliquelemma}, denoted by $\alpha_i\ (1\leq i\leq 3t-4)$.
It is not difficult to derive that both $S_2\xi_1$ and $S_2\xi_2$ are orthogonal to $\alpha_i\ (1\leq i\leq 3t-4)$.
Thus, the dimension of the eigenspace for the $\mathcal{L}$-eigenvalue $\frac{1+t}{t}$ is at least $3t-2$, which indicates that the multiplicity of the $\mathcal{L}$-eigenvalue $\frac{1+t}{t}$ is at least $3t-2$. Noting that $3t+1=n$, we finally obtain that the spectrum of $\mathcal{L}(G_2)$ is $\{0^1,\ {(\frac{3}{2(n-1)})}^1,\ {(\frac{9}{2(n-1)})}^1,\  {(\frac{n+2}{n-1})}^{n-3}\}$.

For graph $G_3$, $\mathcal{L}(G_3)$ has an equitable partition with respect to the vertex partition
$$V(G_3)=V(K_t)\cup V(K_{t})\cup V(K_t)\cup \{x\}\cup \{y\}\cup \{z\}.$$
Denote the corresponding characteristic matrix and quotient matrix of $\mathcal{L}(G_3)$ by $S_3$ and $Q_3$ respectively, then
\vspace{-0.6cm}
$$S_3=\bordermatrix{
    & & & & & & \cr
1      & 1        & 0       & 0       & 0      & 0       & 0   \cr
       & \vdots   & \vdots  & \vdots  & \vdots & \vdots  & \vdots\cr
t      & 1        & 0       & 0       & 0      & 0       & 0 \cr
(t+1)  & 0        & 1       & 0       & 0      & 0       & 0 \cr
       & \vdots   & \vdots  & \vdots  & \vdots & \vdots  & \vdots\cr
2t     & 0        & 1       & 0       & 0      & 0       & 0\cr
(2t+1) & 0        & 0       & 1       & 0      & 0       & 0\cr
       & \vdots   & \vdots  & \vdots  & \vdots & \vdots  & \vdots\cr
3t     & 0        & 0       & 1       & 0      & 0       & 0\cr
(3t+1) & 0        & 0       & 0       & 1      & 0       & 0\cr
(3t+2) & 0        & 0       & 0       & 0      & 1       & 0\cr
(3t+3) & 0        & 0       & 0       & 0      & 0       & 1\cr}
,$$
$$Q_3=
\left(
  \begin{array}{cccccc}
    \frac{2}{t+1} & 0 & 0 & \frac{-1}{(t+1)\sqrt{2}} & \frac{-1}{(t+1)\sqrt{2}} & 0 \\
    0 & \frac{2}{t+1} & 0 & 0 & \frac{-1}{(t+1)\sqrt{2}} & \frac{-1}{(t+1)\sqrt{2}} \\
    0 & 0 & \frac{2}{t+1} & \frac{-1}{(t+1)\sqrt{2}} & 0 & \frac{-1}{(t+1)\sqrt{2}} \\
    \frac{-t}{(t+1)\sqrt{2}} & 0 & \frac{-t}{(t+1)\sqrt{2}} & 1 & \frac{-1}{2(t+1)} & \frac{-1}{2(t+1)}\\
    \frac{-t}{(t+1)\sqrt{2}} & \frac{-t}{(t+1)\sqrt{2}} & 0 & \frac{-1}{2(t+1)} & 1 & \frac{-1}{2(t+1)} \\
    0 & \frac{-t}{(t+1)\sqrt{2}} & \frac{-t}{(t+1)\sqrt{2}} & \frac{-1}{2(t+1)} & \frac{-1}{2(t+1)} & 1 \\
  \end{array}
\right).$$
After calculating, the spectrum of $Q_3$ is
$$\{0^1,\ {(\frac{3}{2(t+1)})}^2,\ {(\frac{t+2}{t+1})}^3\},$$
and there are three linearly independent eigenvectors, denoted by $\xi_1, \xi_2$ and $\xi_3$, for the eigenvalue $\frac{t+2}{t+1}$.
By Lemma \ref{quotientmatrixlemma}, all the eigenvalues of $Q_3$  are also the eigenvalues of $\mathcal{L}(G_3)$, and $S_3\xi_1,\ S_3\xi_2$ and $S_3\xi_3$ are linearly independent eigenvectors for the $\mathcal{L}$-eigenvalue $\frac{t+2}{t+1}$.
Moreover, by Lemma \ref{cliquelemma}, $\frac{t+2}{t+1}$ is an $\mathcal{L}$-eigenvalue of $G_3$ with multiplicity at least $3t-3$,  and one can easily obtain the corresponding $3t-3$ eigenvectors for $\frac{t+2}{t+1}$ similar as stated in Lemma \ref{cliquelemma}, denoted by $\alpha_i\ (1\leq i\leq 3t-3)$.
It is clear that $S_3\xi_i\ (1\leq i\leq 3)$ is orthogonal to $\alpha_i\ (1\leq i\leq 3t-3)$.
Hence, the dimension of the eigenspace for the $\mathcal{L}$-eigenvalue $\frac{t+2}{t+1}$ is at least $3t$, which indicates that $m(\frac{t+2}{t+1})\geq 3t$. Since $3t+3=n$ in $G_3$, then the spectrum of $\mathcal{L}(G_3)$ is $\{0^1,\ {(\frac{9}{2n})}^2,\ {(\frac{n+3}{n})}^{n-3}\}$.
\hfill$\square$

\section{Proof of Theorem 1.1}

\quad The proof of Theorem 1.1 will be completed by the following two theorems.

\begin{theorem}\label{theorem1} \ Let $G$ be a connected graph of order $n\geq 5$. Then $G\in \mathcal{G}(n, n-3)$ with $\rho_{n-1}(G)=1$ if and only if $G$ is a complete tripartite graph $K_{a,b,c}$ ($a+b+c=n$) or $K_n-e$.
\end{theorem}

\noindent
{\bf Proof.} \ Firstly, we show the sufficiency part.

If $G=K_{a,b,c}$ with $a+b+c=n$, then by Lemma \ref{twinpointslemma} the multiplicity of $1$ as an $\mathcal{L}$-eigenvalue is at least $a-1+b-1+c-1=n-3$. In addition, Lemma \ref{secondlargesteigenvalue} indicates that $\rho_{2}(G)>1$. Hence, $K_{a,b,c}$ contains $1$ as an $\mathcal{L}$-eigenvalue with multiplicity $n-3$ and $\rho_{n-1}(K_{a,b,c})=1$.

If $G=K_{n}-e$, by applying Lemmas \ref{twinpointslemma} and \ref{cliquelemma}, we get that $1$ and $1+\frac{1}{n-1}$ are $\mathcal{L}$-eigenvalues of $K_n-e$ with multiplicity at least $1$ and $n-3$, respectively. Then the remaining nonzero $\mathcal{L}$-eigenvalue of $K_n-e$ is $n-1-(n-3)(1+\frac{1}{n-1})=\frac{n+1}{n-1}$. Hence, $K_n-e$ contains $1+\frac{1}{n-1}$ as an $\mathcal{L}$-eigenvalue with multiplicity $n-3$ and $\rho_{n-1}(K_{n}-e)=1$.

Next, we present the necessity part.

Suppose that $G\in \mathcal{G}(n, n-3)$ with $\rho_{n-1}(G)=1$,  then $G$ is a complete multipartite graph from Lemma \ref{secondleasteigenvalue}. We first show that $G$ is neither a complete graph nor a complete bipartite graph. Since $G\in \mathcal{G}(n, n-3)$ is a connected graph, then $G$ has at least three distinct $\mathcal{L}$-eigenvalues. Thus $G$ is not a complete graph (containing two distinct $\mathcal{L}$-eigenvalues). Further, if $G$ is a complete bipartite graph, then $\rho_{2}(G)=1$ from Lemma \ref{secondlargesteigenvalue}, which implies that $1$ is an $\mathcal{L}$-eigenvalue of multiplicity $n-2$, and thus $G\notin \mathcal{G}(n, n-3)$, a contradiction. Next, the remaining proof can be divided into the following two cases.

\vskip 2mm
\noindent
{\bf Case 1.} \ Suppose that the multiplicity of $1$ as an $\mathcal{L}$-eigenvalue is $n-3$, i.e., $m(1)=n-3$.

In this case, it is obvious that
$$\rho_{1}(G)\geq \rho_{2}(G)>\rho_{3}(G)=\rho_{4}(G)\cdots =\rho_{n-1}(G)=1>\rho_{n}(G)=0.$$
Recalling that $G$ is a complete multipartite graph, we first  assume that $G$ is a complete $k$-partite graph with $4\leq k\leq n-1$. Choosing four vertices, say $\{v_1, v_2, v_3, v_4\}$, from four distinct partite, we write the principal submatrix of $\mathcal{L}(G)$ indexed by $v_i\ (1\leq i\leq 4)$ as
$$M=\left(
       \begin{array}{cccc}
         1 & \frac{-1}{\sqrt{d_{v_1}d_{v_2}}} & \frac{-1}{\sqrt{d_{v_1}d_{v_3}}} & \frac{-1}{\sqrt{d_{v_1}d_{v_4}}} \\
         \frac{-1}{\sqrt{d_{v_2}d_{v_1}}} & 1 & \frac{-1}{\sqrt{d_{v_2}d_{v_3}}} & \frac{-1}{\sqrt{d_{v_2}d_{v_4}}} \\
         \frac{-1}{\sqrt{d_{v_3}d_{v_1}}} & \frac{-1}{\sqrt{d_{v_3}d_{v_2}}} & 1 & \frac{-1}{\sqrt{d_{v_3}d_{v_4}}} \\
         \frac{-1}{\sqrt{d_{v_4}d_{v_1}}} & \frac{-1}{\sqrt{d_{v_4}d_{v_2}}} & \frac{-1}{\sqrt{d_{v_4}d_{v_3}}} & 1 \\
       \end{array}
     \right).$$
From Lemma \ref{Interlacinglemma}, we have $$1=\rho_{3}(G)\geq \lambda_3(M)\geq \rho_{n-1}(G)=1,$$
then the third largest eigenvalue $\lambda_3(M)$ of $M$ is equal to 1. As a result, the matrix $M-I_4$ has 0 as an eigenvalue, where $I_4$ is the identity matrix of order $4$. However, one can easily obtain that $det(M-I_4)=\frac{-3}{d_{v_1}d_{v_2}d_{v_3}d_{v_4}}\neq 0$, a contradiction. Therefore, the complete $k$-partite graphs with $4\leq k\leq n-1$ do not belong to $\mathcal{G}(n, n-3)$, which yields that $G$ must be a complete tripartite graph for this case.

\vskip 2mm
\noindent
{\bf Case 2.} \ The multiplicity of $1$ as an $\mathcal{L}$-eigenvalue is not $n-3$, i.e., $m(1)\neq n-3$.

In this case, it is easy to know that $1\leq m(1)\leq 2$ and $G$ is not a complete tripartite graph from the above discussion. Now, let $G$ be a complete $k$-partite graph with $4\leq k\leq n-1$. From Lemmas \ref{twinpointslemma} and \ref{cliquelemma}, we see that there exist at most two partite of $G$ containing more than one vertex, and each of such partite contains at most three vertices.

First, let $N_1$ and $N_2$ be two partite of $G$ such that  $|N_1|\geq 2$ and $|N_2|\geq 2$. Then it follows from Lemma \ref{twinpointslemma} and $1\leq m(1)\leq 2$ that $|N_1|=|N_2|=2$ and each of the remaining $k-2$ partite contains one vertex. Thus $m(1)=2$, and from Lemma \ref{cliquelemma}, $G$ contains $1+\frac{1}{n-1}$ as an $\mathcal{L}$-eigenvalue for $k-2\geq 2$. Therefore, $m(1+\frac{1}{n})=n-3$. By considering the trace of $\mathcal{L}(G)$, we derive that $n=2+(n-3)(1+\frac{1}{n-1})$, which cannot hold.

Second, suppose there exists precisely one partita, say $N_1$, of $G$ containing more than one vertex. Since $1\leq m(1)\leq 2$, then $|N_1|\leq 3$ from Lemma \ref{twinpointslemma}. If $|N_1|=3$, then $m(1)=2$ and one can obtain a contradiction with similar discussion above. As a result, $|N_1|=2$, that is, $G$ is the graph $K_n-e$.

Consequently, the proof is completed.\hfill$\square$

\begin{theorem}\label{theorem2}\ Let $G$ be a connected graph of order $n\geq 5$. Then $G\in \mathcal{G}(n, n-3)$ with $\rho_{n-1}(G)\neq 1$ and $\nu(G)\neq 2$ if and only if $G\in \{G_1, G_2, G_3\}$ (see Fig. 1).
\end{theorem}

\noindent
{\bf Proof.} \ We first present the sufficiency part.
For graphs $G_i\ (1\leq i \leq 3)$ in Fig. 1, since they are not complete multipartite graphs obviously, then $\rho_{n-1}(G_i)\neq 1\ (1\leq i \leq 3)$ by Lemma \ref{secondleasteigenvalue}. In addition, $\nu(G_i)=3\neq 2\ (1\leq i \leq 3)$ clearly. From Lemma \ref{spectrum}, we see that $G_i\in \mathcal{G}(n, n-3)\ (1\leq i \leq 3)$.

Next, we demonstrate the necessity part.
Let $\theta$ be the $\mathcal{L}$-eigenvalue of $G$ with multiplicity $n-3$. Since $G\in \mathcal{G}(n, n-3)$ and $\rho_{n-1}(G)\neq 1$, then $\theta\neq 1$ and $\nu(G)\leq 3$ from Lemma \ref{less3}. Moreover, if $\nu(G)= 1$, then $G$ is the complete graph $K_n$, which contains precisely two distinct  $\mathcal{L}$-eigenvalues ($0$ and $\frac{n}{n-1}$). Therefore, we only need to consider the case of $\nu(G)=3$.
Now, suppose that $\nu(G)=3$ and $\{u, v, w\}$ is a maximum independent set of $G$. Then any vertex out of $\{u, v, w\}$ must be adjacent to at least one of $\{u, v, w\}$. The remaining proof can be divided into the following two cases.

\vskip 2mm
\noindent
{\bf Case 1.} \ Suppose there exists a vertex, say $z$, such that $z$ is adjacent to each of $\{u, v, w\}$.

If this is the case, then applying Lemma \ref{commonvertex} (i) and (iv) we obtain that there exists no vertex adjacent to two of $\{u, v, w\}$ and $z$ is the only vertex adjacent to each of $\{u, v, w\}$. As a result, by Lemma \ref{commonvertex} (ii), $G$ is isomorphic to $\Gamma_1$ in Fig. 2. Next, we show that $s=t=p$ in $\Gamma_1$, i.e., $G$ is isomorphic to $G_1$ in Fig. 1.

Suppose without loss of generality that $s\geq t\geq p$ and $s\geq 2$ in $\Gamma_1$. Then from Lemma \ref{cliquelemma},  $1+\frac{1}{d_u}$ is an $\mathcal{L}$-eigenvalue of $G$.
Denote by $M$ the principal submatrix of $\mathcal{L}(G)-\theta I$ indexed by $\{u, v, w, z\}$, then
$$M=
\left(
  \begin{array}{cccc}
    1-\theta & 0 & 0 & \frac{-1}{\sqrt{d_ud_z}} \\
    0 & 1-\theta & 0 & \frac{-1}{\sqrt{d_vd_z}} \\
    0 & 0 & 1-\theta & \frac{-1}{\sqrt{d_wd_z}} \\
     \frac{-1}{\sqrt{d_zd_u}}& \frac{-1}{\sqrt{d_zd_v}} & \frac{-1}{\sqrt{d_zd_w}} & 1-\theta \\
  \end{array}
\right).
$$
Similarly, we can let $R_z=aR_u+bR_v+cR_w$. Then applying this equation to all the columns of $M$, one can derive that
\begin{equation}\label{e9}
  \begin{array}{rcl}
  (1-\theta)^2 &=& \frac{1}{d_z}(\frac{1}{d_u}+\frac{1}{d_v}+\frac{1}{d_w})\\
   &=& \frac{1}{n-1}(\frac{1}{d_u}+\frac{1}{d_v}+\frac{1}{d_w}).

    \end{array}
\end{equation}

Now, first assume that $\theta=1+\frac{1}{d_u}$. Then it follows from Eq. (\ref{e9}) and the fact $d_u+d_v+d_w=n-1$ that
\begin{equation}\label{e10}
  d_u^2=d_vd_w.
\end{equation}
Recalling that $s\geq t\geq p$ (i.e., $d_u\geq d_v\geq d_w$), we get $d_u=d_v=d_w$ (i.e., $s= t= p$) by Eq. (\ref{e10}), that is, $G$ is isomorphic to $G_1$ in Fig. 1.

Second, assume that $\theta\neq 1+\frac{1}{d_u}$. Then from Lemma \ref{cliquelemma},  $m(1+\frac{1}{d_u})\geq s-1$. Note that the least $\mathcal{L}$-eigenvalue of $G$ is $0$. Thus $s-1\leq m(1+\frac{1}{d_u})\leq 2$, which implies $s\leq 3$. Moreover, if $s=3$, then $t<3$; otherwise, $s=t=3$, then $d_u=d_v$ and  $m(1+\frac{1}{d_u})\geq 4$ by Lemma \ref{cliquelemma}, a contradiction. As a result, we only need to check the graphs in the five cases
\begin{equation*}
\begin{cases}
  \{s=3, t=p=2\} \\
  \{s=3, t=2, p=1\}\\
  \{s=3, t=p=1\}\\
  \{s=t=2, p=1\}\\
  \{s=2, t=p=1\},
\end{cases}
\end{equation*}
none of which belong to $\mathcal{G}(n, n-3)$ by a direct calculation.

\vskip 2mm
\noindent
{\bf Case 2.} \ Suppose there exists no vertex adjacent to each of $\{u, v, w\}$.

In this case, we can see that the vertices out of $\{u, v, w\}$ are adjacent to precisely either one of $\{u, v, w\}$ or two of $\{u, v, w\}$. Since $G$ is connected, by Lemma \ref{commonvertex} (iii), there are two subcases to be considered.

\vskip 2mm
\noindent
{\bf Subcase 2.1} \ There is precisely one pair of $\{u, v, w\}$ having no common vertex.

Suppose without loss of generality that $u$ and $w$ have no common vertex. In other words, there exist $x$ and $y$ such that $x\thicksim u$, $x\thicksim v$, $x\nsim w$ and $y\thicksim v$, $y\thicksim w$, $y\nsim u$. Then by Lemma \ref{commonvertex} (v), we obtain $x\thicksim y$ and
\begin{equation}\label{e11}
d^2_v=d_ud_x-d_ud_v=d_yd_w-d_wd_v.
\end{equation}
Further, it follows from Lemma \ref{commonvertex} (i) that $x$ (resp., $y$) is the unique vertex adjacent to $u, v$ (resp., $v, w$). Therefore, $G$ is isomorphic to $\Gamma_2$ (see Fig. 2) from  Lemma \ref{commonvertex} (ii). Focusing on graph $\Gamma_2$, one can easily obtain
\begin{equation}\label{e12}
\begin{cases}
d_x=p+s+1\\ d_y=s+t+1\\ d_u=p\\ d_v=s+1\\ d_w=t,
\end{cases}
\end{equation}
which indicate that
\begin{equation}\label{e13}
\begin{cases}
  d_x = d_u+d_v \\
  d_y = d_v+d_w.
\end{cases}
\end{equation}
Combining Eqs. (\ref{e11}) and (\ref{e13}), we derive that $d_u=d_v=d_w$, which together with Eq. (\ref{e12}) yields that
$p=s+1=t$. Consequently, $G$ is isomorphic to $G_2$ in Fig. 1.

\vskip 2mm
\noindent
{\bf Subcase 2.1} \ Any two of $\{u, v, w\}$ have a common vertex.

For this case, let $S_{\{u,v\}}$ (resp., $S_{\{v,w\}}$ and $S_{\{u,w\}}$) be the set of the common vertices of $\{u,v\}$ (resp., $\{v,w\}$ and $\{u,w\}$).  Lemma \ref{commonvertex} (i) tells us that $$|S_{\{u,v\}}|=|S_{\{v,w\}}|=|S_{\{u,w\}}|=1,$$
and then let $x\in S_{\{u,v\}}$, $y\in S_{\{v,w\}}$ and $z\in S_{\{u,w\}}$. For any two of $\{x, y, z\}$, Lemma \ref{commonvertex} (v) indicates that $x\thicksim y$, $y\thicksim z$, $x\thicksim z$ and
\begin{equation}\label{e14}
  d^2_v=d_ud_x-d_ud_v=d_yd_w-d_wd_v.
\end{equation}
Then now each of the vertices out of $\{u, v, w, x, y, z\}$ is adjacent to precisely one of $\{u, v, w\}$. Thus, applying Lemma \ref{commonvertex} (ii), we obtain that $G$ is isomorphic to $\Gamma_3$ (see Fig. 2). In graph $\Gamma_3$, one can easily derive that
\begin{equation*}
\begin{cases}
d_x=d_u+d_v\\
d_y=d_v+d_w,
\end{cases}
\end{equation*}
which together with Eq. (\ref{e14}) yields that $d_u=d_v=d_w$ (i.e., $p=s=t$). Therefore, $G$ is isomorphic to $G_3$ in Fig. 1.

As a consequence, the proof is finished. \hfill$\square$

\begin{figure}[htbp]
  \centering
  \setlength{\abovecaptionskip}{0cm} 
  \setlength{\belowcaptionskip}{0pt}
  \includegraphics[width=5 in]{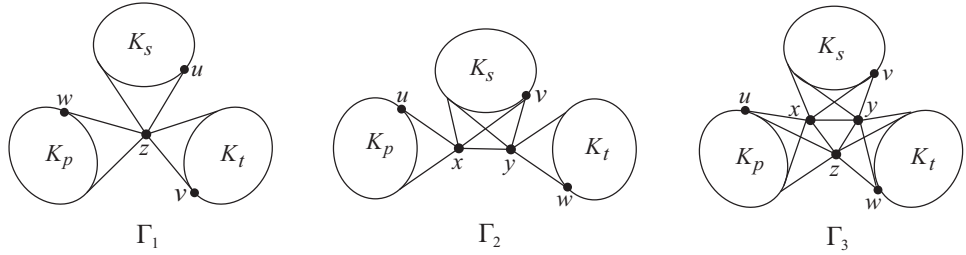}
  \caption{The graphs $\Gamma_1\ (d_z=n-1)$, $\Gamma_2$ and $\Gamma_3$.}
\end{figure}

From Theorem \ref{theorem1} and Lemma \ref{spectrum}, the spectra of the graphs $K_{a,b,c}$, $K_n-e$, and $G_i$ $(1\leq i\leq 3)$ in Fig. 1 are distinct. Then the following corollary is clear.

\begin{corol}\ Suppose that $G\in \mathcal{G}(n, n-3)$ with  $\rho_{n-1}(G)=1$ or $\rho_{n-1}(G)\neq 1$ and $\nu(G)\neq 2$, then $G$ is determined by its spectrum.
\end{corol}

\begin{remark}\ To characterize all the graphs of $\mathcal{G}(n, n-3)$ completely, there is only one remaining case to be considered:  $\rho_{n-1}(G)\neq 1$ and $\nu(G)= 2$.
For this case, we find the following two graphs in Fig. 3 and conjecture that there exists no other graphs belonging to $\mathcal{G}(n, n-3)$.
\end{remark}

\begin{figure}[htbp]
  \centering
  \setlength{\abovecaptionskip}{0cm} 
  \setlength{\belowcaptionskip}{0pt}
  \includegraphics[width=2.8 in]{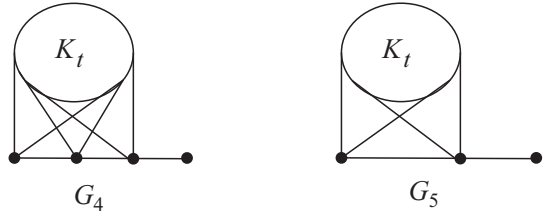}
  \caption{The graphs $G_4$ and $G_5$.}
\end{figure}

\end{document}